\documentclass[12pt]{amsart}

\newtheorem{defin}{DEFINITION}
\newtheorem{prop}{PROPOSITION}
\newtheorem{lemma}{LEMMA}
\newtheorem{cor}{COROLLARY}

\newcommand{\cal}{\mathcal}
\newcommand{\IN}{I\mskip-6mu N}
\newcommand{\IR}{I\mskip-6mu R}
\newcommand{\e}{{\cal E}}

\newcommand{\df}{(\e,D(\e))}
\newcommand{\dif}{{\mbox{$X_t$}}}
\newcommand{\po}{{\pi_\sigma}}
\newcommand{\gd}{{\nabla^\Gamma}}
\newcommand{\fcb}{{\cal F}C^\infty_b}

\begin{document}
\title[A support property for interacting diffusion processes]  
{A support property for infinite dimensional interacting diffusion
  processes.  
\vskip .4cm\noindent
 Une propri\'et\'e de support pour des processus de diffusion  
en dimension infinie avec interaction.}
\author{Michael R\"ockner} 
\address{\hskip-\parindent
 Michael R\"ockner\\
 Fakult\"at f\"ur Mathematik\\
 Universit\"at Bielefeld\\
 Postfach 10 01 31\\
 D-33501 Bielefeld, Germany}
 \email{roeckner@mathematik.uni-bielefeld.de}

\author{Byron Schmuland} 
 \address{\hskip-\parindent
Byron Schmuland\\
 Department of Mathematical Sciences\\
 University of Alberta\\
 Edmonton, Alberta, Canada
 T6G 2G1}
 \email{schmu@stat.ualberta.ca}

\thanks{Research at MSRI is supported in part by NSF grant DMS-9701755.}

\begin{abstract}
The Dirichlet form associated with the intrinsic gradient on Poisson space is
known to be quasi-regular on the complete metric space $\ddot\Gamma=$
$\{Z_+$-valued  
Radon measures on $\IR^d\}$. We show that under mild conditions, the set $\ddot\Gamma\setminus\Gamma$
 is $\e$-exceptional, where $\Gamma$ is the space of locally finite configurations in $\IR^d$, that is,
measures $\gamma\in\ddot\Gamma$ satisfying $\sup_{x\in\IR^d}\gamma(\{x\})\leq 1$.
Thus, the associated diffusion lives on the smaller space $\Gamma$. This result also holds 
for Gibbs measures with superstable interactions.

\vskip .3cm\noindent
{\sc R\'esum\'e.}  Il est connu que la forme de Dirichlet
  associ\'ee au gradient intrins\`eque sur l'espace de Poisson 
 est quasi-r\'eguli\`ere sur l'espace m\'etrique complet 
$\ddot\Gamma=\{$ mesures de Radon sur $\IR^d$ \`a valeurs 
dans $Z_+\}$. Sous des conditions faibles,
on montre que
 l'ensemble $\ddot\Gamma\setminus\Gamma$
  est $\e$-exceptionnel, o\`u $\Gamma$ d\'esigne 
  l'ensemble de toutes les configurations localement finies dans $\IR^d$,
 c'est-\`a-dire, les mesures
  $\gamma\in\ddot\Gamma$ qui verifient $\sup_{x\in\IR^d}\gamma(\{x\})\leq 1$.
  La diffusion associ\'ee 
 prend donc ses valeurs dans l'espace des phases $\Gamma$.
 Ce r\'esultat est valable meme si 
 $\mu$ est une mesure de Gibbs associ\'ee \`a un potentiel
 superstable.
\end{abstract}

\maketitle

\vskip .6cm\noindent
{\bf Version fran\c{c}aise abr\'eg\'ee}
\vskip .3cm
 Soit $\Gamma$ l'ensemble de toutes les configurations localement finies dans $\IR^d$.
$\Gamma$ est vu comme un sous-ensemble de l'espace des mesures de Radon sur 
$\IR^d$ muni de la topologie vague.
Soit $\sigma(dx)=\rho(x)m(dx)$, $m$ \'etant la mesure de Lebegue sur $\IR^d$ et
$\rho^{1/2}\in H^{1,2}_{\rm loc}(\IR^d)$, $\rho(x)>0$ $m$-presque tout $x\in \IR^d$.
Soit $\pi_\sigma$ la mesure de Poisson d'intensit\'e $\sigma$ sur $\Gamma$.
On appellera mesure de Poisson mixte toute mesure $\mu$ de la forme
 $\mu := \int_{\IR_+} \pi_{z\sigma}\,\lambda(dz).$
D\'efinissons l'espace des fonctions cylindriques lisses et born\'ees  
\[\fcb := \left\{  u \mid u(\gamma) =
g(\langle f_1,\gamma\rangle, \langle f_2,\gamma\rangle, \dots, \langle f_n,\gamma\rangle)\right\},\] 
o\`u $f_i\in C^\infty_0(\IR^d)$, $g\in C^\infty_b(\IR^n)$, et
 $\langle f,\gamma\rangle = \sum_{x\in\gamma} f(x)$.
D'apr\`es Albeverio, Kondratiev, et R\"ockner \cite{1} nous d\'efinissons le gradient intrins\`eque d'une
fonction $u\in\fcb$ par
\[\left(\gd u\right)(\gamma ; x):= \sum_{i=1}^n \frac{\partial g}{\partial x_i}
\left(\langle f_1,\gamma\rangle, \langle f_2,\gamma\rangle, \dots, \langle f_n,\gamma\rangle\right)
\nabla f_i(x).\]
Donc on a l'op\'erateur carr\'e du champ 
 \[\Gamma(u,v)(\gamma) 
 := \int_{\IR^d} \langle (\gd u)(\gamma;x), (\gd v)(\gamma;x)\rangle_{\IR^d} \,\gamma(dx).\]
La forme de Dirichlet engrendre\'e par cet op\'erateur carr\'e du champ est la fermeture $\df$ de la forme
$(\e,\fcb)$ d\'efinie pour $u,v\in\fcb$ par 
\[\e(u,v):= \int_\Gamma \Gamma(u,v)(\gamma)\,\mu(d\gamma).\]
 
On sait que cette forme est quasi-r\'eguli\`ere sur l'espace des phases complet 
$\ddot\Gamma (\supset \Gamma)$. On va prouver ici que 
$\df$ est quasi-r\'eguli\`ere sur $\Gamma$; pour cela il est suffisant de voir
que l'ensemble $\ddot\Gamma\setminus \Gamma$ est $\e$-exceptionnel. 
On a le r\'esultat suivant:
\vskip .2cm
 {PROPOSITION 1. \it Si $d\geq 2$, $\rho \in L^2_{\rm loc}(dx)$, et $\int_{\IR_+} z^2\,\lambda(dz)<\infty$,
 l'ensemble $\ddot\Gamma\setminus\Gamma$ est $\e$-exceptionnel.}
\vskip .2cm
Ce r\'esultat est valable meme si $\mu$ est une mesure de Gibbs associ\'e \`a un potentiel 
superstable au sens de Ruelle, puisque l'existence des points avec masse $\geq 2$ est un \'ev\'enement 
local et une telle mesure de Gibbs est localement absolument continue par rapport \`a une mesure de Poisson. 

\vskip .3in
\hrule
\vskip .3in
Let $\Gamma$ be the space of locally finite configurations in $\IR^d$,
\begin{equation}
 \Gamma:=\{\gamma \subset \IR^d : |\gamma\cap K|<\infty\hbox{ for every compact }K\}.
\end{equation}
A configuration $\gamma$ will
be identified with the Radon measure 
$ \sum_{x\in \gamma} \varepsilon_{x}$.
 The space $\Gamma$ will be given the topology of vague convergence of measures, and
 measures on $\Gamma$ are defined on the corresponding Borel sets.
Define
\begin{eqnarray}
\nonumber
\fcb := \left\{ \right. u &\mid& \gamma\mapsto u(\gamma) =
g(\langle f_1,\gamma\rangle, \langle f_2,\gamma\rangle, \dots, \langle f_n,\gamma\rangle) \\
 & & \left.\hbox{ for some } f_i\in C^\infty_0(\IR^d)\hbox{ and } g\in C^\infty_b(\IR^n)\right\},
\end{eqnarray}
where  $\langle f,\gamma\rangle = \sum_{x\in\gamma} f(x)$.

 For $u\in\fcb$, we define the gradient $\gd u$ at the point $\gamma\in\Gamma$ as an element of the
``tangent space" $T_\gamma(\Gamma):=L^2(\IR^d\to \IR^d ;\gamma)$ by
\begin{equation}
x\mapsto \left(\gd u\right)(\gamma ; x):= \sum_{i=1}^n \frac{\partial
  g}{\partial x_i} 
\left(\langle f_1,\gamma\rangle, \langle f_2,\gamma\rangle, \dots, \langle f_n,\gamma\rangle\right)
\nabla f_i(x).
\end{equation}
Here $\nabla $ refers to the usual gradient on $\IR^d$.
 It is not hard to prove that $\gd u$ is well-defined, even though the representation of 
$u$ as a cylinder function is not unique.

\begin{defin}{
 For $u,v\in\fcb$ define the {\it square field operator} $\Gamma$ as 
 the real-valued function on $\Gamma$ given by
 \begin{eqnarray}\nonumber
 \gamma\mapsto\Gamma(u,v)(\gamma) 
 &:=&\langle \gd u,\gd v\rangle_{T_\gamma(\Gamma)}\\
 &=& \int_{\IR^d} \langle (\gd u)(\gamma;x), (\gd v)(\gamma;x)\rangle_{\IR^d} \,\gamma(dx).
 \end{eqnarray}
 We will often use the abbreviation $\Gamma(u):=\Gamma(u,u)$.}
\end{defin}

Let $\sigma$ be a measure on $\IR^d$ that has a density $\rho$ with respect to Lebesgue measure
 satisfying $\rho >0$ almost everywhere, and $\rho^{1/2}\in H^{1,2}_{loc}(\IR^d)$. Here
 $H^{1,2}_{loc}(\IR^d)$ denotes the local Sobolev space of order 1 in $L^2_{loc}(\IR^d; dx)$.
 The {\it Poisson measure} $\po$
 with intensity measure $\sigma$ is the probability measure on $\Gamma$ characterized by
 \begin{equation}
\int_\Gamma \exp(\langle f,\gamma\rangle)\,\po(d\gamma)
= \exp\left(\int_{\IR^d} (e^{f(x)}-1)\,\sigma(dx)\right),
\label{poisson}
\end{equation}
 for $f\in C_0(\IR^d)$.
  A {\it mixed Poisson measure} is given by 
\begin{equation}
 \mu := \int_{\IR_+} \pi_{z\sigma}\,\lambda(dz),\label{mix}
\end{equation}
 where $\lambda$ is a probability measure such that
$\int_{\IR_+}z\,\lambda(dz)<\infty$.
\begin{defin}
 For $u,v\in\fcb$ we define the pre-Dirichlet form 
\begin{equation} \e(u,v):= \int_\Gamma \Gamma(u,v)(\gamma)\,\mu(d\gamma).\label{form} \end{equation}
\end{defin}

Albeverio, Kondratiev, and R\"ockner proved an integration by parts formula
 \cite[Corollary 4.1 and Remark 4.3]{3} which implies that the form $(\e,\fcb)$ is closable, 
and that its closure $\df$ is a symmetric, local, Dirichlet form.
The quasi-regularity of $\df$ has been proven for certain cases 
 by Yoshida \cite{9}, and in general by Ma and R\"ockner \cite{4}. 
However, $\Gamma$ is not completely metrizable with respect to the
vague topology so it is necessary to use the Polish state space
\begin{equation}
 \ddot\Gamma := \{ Z_+\hbox{-valued Radon measures on }\IR^d\}.
\end{equation}
Since $\Gamma\subset \ddot\Gamma$ and ${\cal B}(\ddot\Gamma)\cap \Gamma = {\cal B}(\Gamma)$, we can consider
$\mu$ as a measure on $(\ddot\Gamma, {\cal B}(\ddot\Gamma))$ and correspondingly $\df$ as a 
Dirichlet form on $L^2(\ddot\Gamma; \mu)$.
The associated  Markov process
$( (\dif)_{t\geq 0}, (P_\gamma)_{\gamma\in\ddot\Gamma})$
 has vaguely continuous sample paths since $\df$ is a local form 
(cf. \cite[Chapter V, Theorem 1.11]{6}).

We will use the following lemma from Dirichlet form theory.
\begin{lemma} {
 Let $u_n \in D(\e)$ be continuous functions,
$\sup_n\e(u_n,u_n)<\infty$, and $u_n\to u$ pointwise.  Then $u$
is an $\e$-quasi-continuous function, in particular,
for $\mu$-almost every $\gamma\in\ddot\Gamma$,
\begin{equation}
  P_\gamma\left( t\to u(\dif)\hbox{ is continuous }\right) =1.\label{continuity}
\end{equation}} 
\end{lemma}

For example, if $1_N$ is $\e$-quasi-continuous and $\mu(N)=0$, then
for $\mu$-almost every $\gamma\in\ddot\Gamma$, 
\begin{equation}
 P_\gamma\left(X_t\not\in N\hbox{ for all } 0<t<\infty\right)=1.
\end{equation}
Such a set $N$ is called {\it $\e$-exceptional}.
\vskip .3cm\noindent
\begin{prop}
 {If  $d\geq 2$, $\rho \in L^2_{\rm loc}(dx)$, and $\int_{\IR_+} z^2\,\lambda(dz)<\infty$,
 then the set $\ddot\Gamma\setminus\Gamma$ is $\e$-exceptional.}\label{main_result}
\end{prop}
\par\par\vskip .1cm
{\it Proof. }
 It suffices to prove the result locally, that is, to show that for every positive
integer $a$, the function $u:=1_N$ is quasi-continuous, where
\begin{equation}
 N:=\{ \gamma : \sup(\gamma(\{x\}): x\in [-a, a]^d) \geq 2\}.\label{N}
\end{equation}

 Our analysis begins with a smooth partition of $\IR^d$ into small pieces.
 Let $\phi$ be a $C^\infty_0(\IR)$ function satisfying $1_{[0,1]}\leq \phi\leq 1_{[-1/2, 3/2)}$ and
 $|\phi^\prime|\leq 3\times 1_{[-1/2,3/2)}$, and for any $n\in\IN$ and $i=(i_1,\dots,i_d)\in Z^d$,
define a $C^\infty_0(\IR^d)$ function by 
\begin{equation}
 \phi_i(x):=\prod^d_{k=1} \phi(nx_k-i_k).
\end{equation}
We also let $I_i(x):=\prod^d_{k=1} 1_{[-1/2,3/2)}(nx_k-i_k)$ and note that $\phi_i\leq I_i$.
 Taking the $j^{\rm th}$  partial derivative of $\phi_i$ gives
\begin{equation}
 \partial_j \phi_i(x)=n\phi^\prime(nx_j-i_j)\prod_{k\not= j} \phi(nx_k-i_k),
\end{equation}
and so $(\partial_j \phi_i(x))^2\leq 9n^2 I_i(x)$. Adding over $j$ from 1 to $d$ gives
us \begin{equation}
|\nabla \phi_i(x)|^2\leq 9 n^2 d I_i(x). \label{bound1}
\end{equation}
Let $\psi$ be a smooth function on $\IR$ satisfying $1_{[2,\infty)}\leq \psi\leq 1_{[1,\infty)}$
 and $|\psi^\prime|\leq 2\times 1_{(1,\infty)}$.
Choosing $A:=  Z^d\cap [-na, na]^d$, 
 define a continuous element of $D(\e)$ by
\begin{equation}
u_n(\gamma):=\psi\bigg(\sup_{i\in A} \langle\phi_i, \gamma\rangle\bigg).\label{un}
\end{equation}
Then $u_n\to u$ pointwise as $n\to \infty$, so  
 to apply Lemma 1 
we only need prove that $\sup_n \e(u_n, u_n) < \infty$.
We begin by bounding $\Gamma(u_n)$, the
square field operator applied to $u_n$. First note that 
\begin{equation}
\left(\psi^\prime\bigg(\sup_{i\in A} \langle\phi_i, \gamma\rangle\bigg)\right)^2 
\leq 4\times 1_{(\sup_{i\in A} \langle\phi_i, \gamma\rangle > 1)}
\leq 4\times 1_{(\sup_{i\in A} \langle I_i, \gamma\rangle \geq 2)},\label{bound2}
\end{equation}
where for the final inequality we use the fact that $\langle I_i,\gamma\rangle$ is an integer.
 Therefore, first  using the inequality 
$\Gamma(u\vee v)\leq \Gamma(u)\vee\Gamma(v)$,
and then using 
  (\ref{bound1}) and (\ref{bound2}), we get 
\begin{eqnarray}
\nonumber\Gamma(u_n)(\gamma) &=& \left(\psi^\prime\bigg(\sup_{i\in A} \langle\phi_i, 
\gamma\rangle\bigg)\right)^2 \Gamma(\sup_{i\in A} \langle\phi_i,\cdot \rangle) (\gamma) \\ \label{gun}
&\leq& \left(\psi^\prime\bigg(\sup_{i\in A} \langle\phi_i, \gamma\rangle\bigg)\right)^2
 \sup_{i\in A} \Gamma(\langle\phi_i,\cdot \rangle) (\gamma)\label{bound3} \\
\nonumber &=& \left(\psi^\prime\bigg(\sup_{i\in A} \langle\phi_i, \gamma\rangle\bigg)\right)^2
\sup_{i\in A} \int|\nabla\phi_i(x)|^2\,\gamma(dx) \\
\nonumber &\leq& 4\times 1_{(\sup_{i\in A}\langle I_i,\gamma\rangle \geq 2)}\, 
9 \,n^2\, d\,\sup_{i\in A}\langle I_i,\gamma\rangle \\
\nonumber &\leq& 36 n^2 d \sum_{i\in A} 1_{(\langle I_i, \gamma\rangle \geq 2)} \langle I_i,\gamma\rangle.
\end{eqnarray}
 From equation (\ref{mix}) we have
 \begin{eqnarray}
 \int_{(\langle I_i,\gamma\rangle\geq 2)} \langle I_i,\gamma\rangle \mu(d\gamma) 
 &=& \int_{\IR_+} z \langle I_i,\sigma\rangle \left(1-e^{-z\langle I_i,\sigma\rangle}\right)\,\lambda(dz)\\
 &\leq& \langle I_i,\sigma\rangle^2 \int_{\IR_+} z^2 \,\lambda(dz),
 \end{eqnarray}
 and combined with (\ref{bound3}) this gives
\begin{equation}
\e(u_n, u_n)\leq c\,n^2 \sum_{i\in A} \langle I_i,\sigma\rangle^2. \label{ebound}
\end{equation}

Although the supports of the indicator functions $I_i$ are not disjoint, 
 each point belongs to at most $2^d$ of the sets $\{I_i=1\}$ for $i\in A$.
Therefore the Cauchy-Schwarz inequality gives us
\begin{eqnarray}
\sum_{i\in A}  \langle I_i,\sigma\rangle^2
 &=& \sum_{i\in A} \left(\int {I_i}(x)\, \rho(x)\,dx\right)^2 \\
\nonumber&\leq& \sum_{i\in A} \left(\int {I_i}(x)\, \rho(x)^2 \,dx\right) \,
          \left(\int I_i(x) \,dx\right)\\
\nonumber&\leq& 2^d \int_{[-(a+1), a+1]^d} \rho(x)^2 \,dx \, (2/n)^d,
\end{eqnarray}
and combining this with (\ref{ebound}) 
we find that 
\begin{equation}
\e(u_n, u_n)\leq c n^{2-d}.
\end{equation}
Since  $d\geq 2$ 
we see that $\sup_n \e(u_n, u_n)< \infty$.

\begin{cor} The conclusion of Proposition 1 holds 
 if $\mu$ is replaced by  
 any Gibbs measure $\nu$, in the sense of Ruelle \cite{8}, 
 with a pair potential $\Phi$ 
that is superstable, lower regular, and
 where $\int |\exp(-\Phi(x))-1|\,dx<\infty$.
\end{cor}
\par\par\vskip .1cm
{\it Proof. }
A Gibbs measure $\nu$  as in the statement admits a system of 
density distributions that satisfy $\sigma_\Delta^n(x_1,\dots,x_n)\leq \xi^n$ 
for some constant $\xi$ \cite[Theorem 5.5]{8}.

Any measure $\nu$ on $\Gamma$ admitting such  density distributions
is locally absolutely continuous with respect to the Poisson measure $\pi_{\xi m}$.
The bounds obtained in the proof of Proposition 1 are valid for this Poisson
measure ($\lambda=\varepsilon_\xi$, $\sigma=m=$ Lebegue measure) and carry over
trivially to the measure $\nu$.
\vskip .3cm\noindent
{\it Remarks.\ \ } 

1. More detailed analysis shows that the condition $\int z^2\,\lambda(dz)<\infty$
 in Proposition 1 can be dropped 
 if the measure $\sigma$ satisfies 
the growth condition
$ \sigma(S_r)\leq a \exp(br)$
for some $a,b>0$,
where $S_r$ is the sphere of radius $r$ centered at the origin.

2. We point out that the form $(\e,{\cal D})$ in Osada \cite{7}
extends our $\df$, but it is not known whether the forms coincide.
A form with a larger domain has more exceptional sets, so proofs 
of exceptionality are easier in Osada's setting, and do not imply
exceptionality in our setting.

3. The results of this paper have been announced in lectures at 
Osaka, Seoul, and Taipei (July 1997), and at MSRI Berkeley (November 1997).
\vskip .3cm\noindent
{\it Acknowledgements.\ \ }   

 The authors gratefully acknowledge the support of 
the DFG through SFB 343 Bielefeld,
 and of NSERC (Canada). 
We also thank the 
Mathematics Sciences Research Institute 
at Berkeley for a very pleasant stay.

\vskip .5cm


\begin{thebibliography}{AKR9}
\bibitem{1} S.\ Albeverio, Yu.\ G.\ Kondratiev, and M.\ R\"ockner:
`Differential geometry of Poisson spaces', {\it 
 Comptes Rendus de L'Acad\'emie des Sciences Paris}, 
{\bf 323} 1129--1134, 1996.

\bibitem{2} S.\ Albeverio, Yu.\ G.\ Kondratiev, and M.\ R\"ockner:
`Canonical Dirichlet operator and distorted Brownian motion on Poisson spaces',
{\it Comptes Rendus de L'Acad\'emie des Sciences Paris}, {\bf 323} 1179--1184, 1996.

\bibitem{3} S.\ Albeverio, Yu.\ G.\ Kondratiev, and M.\ R\"ockner:
`Analysis and geometry on configuration spaces', 
to appear in the {\it Journal of Functional Analysis}.

\bibitem{4} S.\ Albeverio, Yu.\ G.\ Kondratiev, and M.\ R\"ockner:
`Analysis and geometry on configuration spaces: The Gibbsian case',
to appear in the {\it Journal of Functional Analysis}.

\bibitem{5}
Z.M.~Ma and M.\ R\"ockner: 
`Construction of diffusions on configuration spaces', Preprint.

\bibitem{6}
Z.M.~Ma and M.\ R\"ockner: {\it Introduction to the Theory of (Non-Symmetric)
Dirichlet Forms.} Berlin: Springer 1992.
 
\bibitem{7} H.~Osada: `Dirichlet form approach to infinite-dimensional
Wiener processes with singular interactions', 
{\it Communications in Mathematical
Physics} {\bf 176} 117--131, 1996.

\bibitem{8} 
D.\ Ruelle: `Superstable interactions in classical statistical mechanics',
{\it Communications in Mathematical Physics} {\bf 18} 127--159, 1970.

\bibitem{9}
M.\ W.\ Yoshida: `Construction of infinite dimensional interacting diffusion 
processes through Dirichlet forms', 
{\it Probability Theory and Related Fields}
{\bf 106} 265--297, 1996.
\end{thebibliography}
\end{document}